\documentclass[11pt,amstex]{article}
\usepackage{amsmath}
\usepackage{amsfonts}
\usepackage{amscd}
\usepackage[dvips]{color}
\usepackage[dvips]{color}

\def\a{\alpha}
\def\b{\beta}
\def\ga{\gamma}

\def\la{\lambda}

\def\phi{\varphi}

\def\w{\omega}

\def\co{{\mathcal O}}

\def\Ka{K_\alpha}

\def\be{\begin{equation}}
\def\ee{\end{equation}}
\def\bear{\begin{eqnarray}}
\def\eear{\end{eqnarray}}
\def\best{\begin{eqnarray*}}
\def\eest{\end{eqnarray*}}

\newtheorem{theorem}{Theorem}[section]
\newtheorem{prop}[theorem]{Proposition}

\newtheorem{lemma}[theorem]{Lemma}
\newtheorem{cor}[theorem]{Corollary}

\newtheorem{defn}[theorem]{Definition}

\newtheorem{Vprin}[theorem]{Vanishing Principle}
\newtheorem{remark}[theorem]{Remark}
\newenvironment{rem}{\begin{remark}\rm}{\end{remark}}
\newtheorem{example}[theorem]{Example}
\newenvironment{ex}{\begin{example}\rm}{\end{example}}

\def\non{\noindent}
\def\pf{\non {\bf Proof. }}
\def\qed{\nopagebreak \hskip .1in { $\Box$ }\penalty10000 \hskip\parfillskip \par  }

\def\ti{\times}
\def\del{\overline \partial}

\def\dim{\mbox{\rm dim }}

\def\ind{\mbox{Index\,}}
\def\ker{\mbox{ker\,}}
\def\cok{\mbox{coker\,}}
\def\sign{\mbox{\rm sign}}

\def\ov#1{\overline{#1}}
\def\Z{{\mathbb Z}}
\def\R{{\mathbb R}}
\def\P{{\mathbb P}}
\def\Q{{\mathbb Q}}
\def\cx{{\mathbb C}}

\def\Ka{K\"{a}hler }

\def\H{{\mathcal H}}

\def\O{{\mathcal O}}

\def\M{{\mathcal M}}
\def\bM{\ov\M}
\def\Mgn{\M_{g,n}}
\def\bMgn{\overline{\M}_{g,n}}

\def\ind{{\rm index\;}}

\newcommand{\CM}{\overline{{\mathcal M}}}

\title{\bf A Structure Theorem for the Gromov-Witten Invariants of
K\"{a}hler Surfaces}
\author{ Junho Lee and Thomas H. Parker\thanks{ partially supported by
the N.S.F.}}
\date{}

\addtocounter{section}{0}
\begin{document}
\maketitle

\begin{abstract}
\medskip

We prove a structure theorem for the Gromov-Witten invariants of
compact K\"ahler surfaces with geometric genus $p_g>0$.  Under the
technical assumption that there is a  canonical divisor that is  a
disjoint union of smooth components, the theorem shows that the GW
invariants are universal functions
determined by the genus  of this  canonical
divisor components and the holomorphic Euler characteristic of  the surface.   We compute special cases of these
universal functions.

\end{abstract}
\vskip.2in


Much of the work on the Gromov-Witten invariants of \Ka surfaces has
focused on rational and ruled surfaces, which have geometric genus
$p_g=0$. This paper focuses on surfaces with $p_g>0$, a  class
that includes most elliptic surfaces and most surfaces
of general type.  In this context we prove a general ``structure
theorem'' that shows (with one technical assumption) how the GW
invariants are completely determined by the local geometry of a
generic canonical divisor.

The structure theorem is a  consequence of a  simple fact:  the
``Image Localization Lemma'' of Section 3.    Given a \Ka surface $X$
and a canonical divisor $D\in |K_X|$,  this lemma shows that the
complex structure $J$ on $X$ can be perturbed to a {\em non-integrable}  almost complex
structure $J_D$ with the property the image of all $J_D$-holomorphic
maps lies in the support of  $D$.   This immediately gives
some striking vanishing theorems  for the GW invariants of  \Ka
surfaces  (see Section 3). More importantly, it implies that
 the Gromov-Witten
invariant of $X$ for genus $g$ and $n$ marked points is a sum
$$
GW_{g,n}(X,A)\ =\ \sum GW^{loc}_{g,n}(D_k, A_k)
$$
over the connected components $D_k$ of $D$ of ``local invariants''
that count the contribution of maps whose image lies in or (after
perturbing to a generic moduli space) near $D_k$.    These local invariants have not been previously defined.  The proof of their existence relies on using non-integrable structures and geometric analysis techniques.

Our structure theorem characterizes  the local invariants by expressing them in terms of usual GW invariants of certain standard surfaces.   For this
we make the mild assumption that one can deform the \Ka structure on
$X$ and choose a  canonical divisor $D$ so that all components of $D$
are smooth. With this assumption,  the restriction of $K_X$ to each
component $D_k$ of multiplicity $m_k$  is the normal bundle $N_k$ to
$D_k$, and
$$
N_k^{m_k+1}=K_{D_k},
$$
that is,  $N_k$ is a holomorphic  root of the
canonical bundle of the curve $D_k$. The
local invariants are given by universal functions
 \bear
\label{IntroDefL} L^i(t)\ \in\  \prod_{g,n} H_*(\bMgn
\times D_k^n)\, [[t, \la]].
\eear
depending only on topological invariants of the pair $(D_k,N_k)$. Four types of such universal functions
are relevant. Each   is defined in terms of the GW invariant of a
rational surface or a local  GW invariant of a line bundle over a
curve (they can also be defined using obstruction bundles, but we
do not pursue that approach here). There is one universal function
$L^0(t)$ for exceptional curves. This enters into the
blow-up formula for GW invariants proved in Section \ref{section5}.
The blow-up formula reduces the computation of the GW invariants
to the case of minimal surfaces.  In light of  the vanishing results given in
Section 3, there are only  two types of minimal surfaces  to
consider: properly elliptic surfaces and surfaces of general type.

A minimal properly elliptic surface can be deformed to guarantee the
existence of a canonical divisor whose support is the union of smooth
elliptic fibers. The structure theorem separates
these into two types\,: regular fibers and multiple fibers with multiplicity
$m\geq 2$; the corresponding universal functions are
$L^1(t)$ and $L^2_m(t)$ respectively.   For a minimal surface of general type all canonical divisors are  connected;  we
assume that  one such divisor $D$ is smooth and reduced.   By the adjunction formula,  $D$ has genus $h=K_X^2+1$.  The GW invariant is then given by one of two universal functions $L^3_{h,\pm}(t)$ for this
$h$.

\smallskip

The  GW invariants of a K\"{a}hler surface $X$ can be regarded as a
power series in formal variables $t_A$ with $A\in H_2(X;\Z)$ as
described in Section 1.  Each smooth component $D_k$ of the canonical
class, we can replace $t$ by $t_{D_k}$ in the appropriate universal
function (\ref{IntroDefL}), taking $h$  to be the genus $D_k$.
Pushing forward under the map
$\left(\iota_{D_k}\right)_* :H_*(\bMgn \times D_k^n)\to H_*(\bMgn
\times X^n)$ induced by the inclusion $D_k\hookrightarrow X$ then gives
a series (suppressing $\left(\iota_{D_k}\right)_* $ in the notation)
$$
L^i(t_{D_k})\ \in\  \prod_{g,n} H_*(\bMgn \times X^n)[[t_{D_k}, \la]]
$$
that gives the local contribution of $D_k$ to the GW series.  For
each surface, only a few of these series are needed.  The
structure theorem lists the possibilities.  The contribution $GW^0_X$
of the class $A=0$ (that is, the evaluation of (\ref{defn.7}) at
$t=0$) must be separated out.   Note that $GW^0_X$ has been
explicitly computed (see
\cite{km}). Also note that surfaces with $p_g>0$ have a unique
minimal model (\cite{BHPV}, pg 243).

\begin{theorem}[Structure Theorem]
\label{introthm}
Let $X$ be a closed  K\"{a}hler surface with $p_g>0$ and smooth canonical
divisor $D$.  Write $D=\sum_i E_i + D'$  where $\{E_i\}$ are the
exceptional curves in $D$.  Then the  GW invariant of $X$ is a sum
$$
GW_X\ =\ GW^0_X\ +\ \sum_{E_i} L^0(t_{E_i})\ +\ GW'_X
$$
where $GW'_X$ is given as follows according to the type of the
minimal model $X'$ of $X$:
\pagebreak
\non

\begin{enumerate}
\item  If $X'$ is K3 or abelian then $GW'_X=0$.
\item   If $\pi:X'\to C$ is properly elliptic, we can assume that
the canonical divisor $D'$ has the form $\sum n_j F_j+\sum (m_k-1)F_k$ for regular fibers $F_j$ and smooth multiple fibers $F_k$
of multiplicity $m_k$. We  then have
$$  GW'_X\  =\ k_{\pi} \,L^1(t_F)\ +\
   \sum_k\,L^2_{m_k}(t_{F_k}) $$
where $F$ is a regular fiber,   $t^{m_k}_{F_k}=t_F$, and $k_\pi=\chi(\O_{X} )-2\chi(\O_C)$.
\item If $X'$ is general type  and we can choose $D'$ to be smooth
with multiplicity 1, then  $D'$ has genus $h=K_{X'}^2+1\geq 2$ and
$$
GW'_{X} =
\begin{cases}
 L^3_{h,+}(t_{D'}) & \mbox{if $\chi(\O_{X})$ is even}\\
  L^3_{h,-}(t_{D'}) & \text{if $\chi(\O_{X})$ is odd.}
  \end{cases}
$$
\end{enumerate}
\end{theorem}
A more detailed version of the structure theorem is given, with
proofs, in Sections~4 -- 7.  Section~8 contains some analytic results
about  the linearization  of the   $J_\a$-holomorphic map equation,
which has some remarkable properties.  Those are used in Sections~9
and 10 to explicitly compute the contribution to the GW invariants of
special types of covers. We do this  for double covers, then for all
etale covers of elliptic fibers.

The structure theorem shows that, under the stated hypotheses,
the GW invariants are determined by the map $H^*(X,\Q)\to H^*(D,\Q)$ induced by the inclusion $D\subset X$ and by the parity of the holomorphic Euler characteristic, which is given in terms of the Betti numbers by  $\chi(\co_{X})=\frac12(1-b^1+b^+)$.   In the case  when $X$ is a simply-connected surface of general type with a smooth reduced canonical divisor, this information is determined by the homology and Seiberg-Witten invariants of $X$, and hence depends only on the differentiable structure of $X$.    Furthermore, the SW invariants are equivalent to  Taubes' $Gr$ invariants, which correspond to a subset of GW invariants \cite{T}, \cite{ip0}.   But we learn from the structure theorem that the full set of GW invariants contain exactly the same information as the $Gr$ and SW invariants.

Our structure theorem applies for surfaces with $p_g>0$.  This  is
exactly the case when the dimension of the spaces of stable maps
differs from the dimension of the generalized Severi variety, and
thus the GW invariants are not enumerative invariants.  The
$J_\a$-holomorphic map equation can also be used to define a set of ``Family
GW invariants''  that are directly related
to
enumerative invariants. That context is explained in \cite{l1},
\cite{l2} and \cite{lp}.

One might hope that the structure theorem  extended to non-K\"{a}hler
symplectic
manifolds with $b^+>1$.   Unfortunately, as M. Usher observed
(\cite{U}, page 4),  McMullen and
Taubes have constructed  a  symplectic four-manifold whose GW
invariant is not  the sum of local invariants supported on the
components of the canonical class.

\medskip

 We thank R. Pandharipande for useful and encouraging conversations,
and in particular for pointing out the role of spin curves.  R. Friedman  generously helped us with elliptic surface theory.
We also thank F. Catanese,  E.
Ionel and Bumsig Kim for helpful comments.

\vskip 1cm
\setcounter{equation}{0}
\section{Gromov-Witten invariants}
\label{section1}
\bigskip

 We will use the
definitions and notation of \cite{ip2} for stable maps and the Gromov-Witten invariants;  these
are based on the approach developed by Ruan-Tian \cite{rt2} and
Li-Tian \cite{LT}.  In summary, the key definitions go as follows. A
{\em bubble domain} $B$ is a finite connected union of smooth
oriented 2-manifolds $B_i$ joined at nodes together with $n$ marked
points, none of which are nodes. Collapsing the unstable components
to points gives a connected domain $st(B)$ with some arithmetic
genus $g$.  Let $\ov{\mathcal U}_{g,n}\to \ov\M_{g,n}$ be the
universal curve over the Deligne-Mumford space of genus $g$ curves
with $n$ marked points.  We can put a complex structure $j$ on $B$
by specifying an orientation-preserving map $\phi_0: st(B) \to
\ov{\mathcal U}_{g,n}$ which is a diffeomorphism onto a fiber of
$\ov{\mathcal U}_{g,n}$.  We will often write $C$ for the curve
$(B,j)$. A $(J,\nu)$-{\em holomorphic map} from $B$ is then a map
$(f,\phi):B \to X\times\ov{\mathcal U}_{g,n}$ where
$\phi=\phi_0\circ st$ and which satisfies
$$
\bar{\partial}_Jf\ =\ \phi^*\nu
$$
(here the  perturbation $\nu$  is a tensor on $X\times\ov{\mathcal
U}_{g,n}$; see \cite{ip2}). Such a  map is a {\em stable map} if the
restriction  of $(f,\phi)$ to each component  of $B$ is non-trivial
in homology. For generic $(J,\nu)$  the moduli space
${\M}_{g,n}(X,A)$ of stable $(J,\nu)$-holomorphic maps representing
a class $A\in H_2(X)$ is a smooth orbifold of (real) dimension \bear
-2K_{X}\cdot A +(\dim X-6)(1-g)+2n. \label{dimM} \eear Its
compactification carries a (virtual) fundamental class whose
pushforward under the map \bear \ov{\M}_{g,n}(X,A) \ \ \overset{{\rm
{st}} \times {\rm{ev}}}\longrightarrow \ \ \ov{\M}_{g,n}\ti X^n
\label{evmap} \eear defined by stabilization and evaluation at the
marked points is the Gromov-Witten invariant
$$
GW_{g,n}(X,A) \in H_*(\ov{\M}_{g,n}\ti X^n).
$$
This is equivalent to the collection of ``GW numbers''
\begin{equation}
\label{GWnumbers}
GW_{g,n}(X,A) (\mu; \gamma^1,\dots, \gamma^n)
\end{equation}
obtained by evaluating on  classes $\mu\in H^*(\bMgn)$ and
$\gamma^j\in H^*(X)$ whose total degree is the dimension (\ref{dimM})
of the space of stable maps.  The number (\ref{GWnumbers})
is obtained by choosing  (generic) geometric representatives
$M\subset\bMgn$ and $\Gamma_i$ of the classes  Poincar\'{e} dual to
$\mu\in H^*(\bMgn)$ and $\ga^i\in H^*(X)$ and counting, with sign,
the finite set of maps $f:C\to X$ in $st(C)\in M$ and  $f(x_i)\in
\Gamma_i$ for   each marked point  $x_i$.

It is convenient to assemble these  into a
single invariant  by  introducing
variables  $\lambda$ to keep track of the Euler class and $t_A$ satisfying
$t_At_{B}=t_{A+B}$ to   keep
track of $A$.  The  GW series  of $(X,\w)$ is then the formal series
\begin{equation}
GW_{X}\ =\ \sum_{A,g,n} \frac{1}{n!}\,GW_{g,n}(X,A) \ \,t_A\ \la^{2g-2}.
\label{defn.7}
\end{equation}
\vskip1cm

\setcounter{equation}{0}
\section{$J_\a$-holomorphic maps into K\"ahler surfaces}
\label{section2}
\bigskip

Fix a K\"ahler surface $(X, J, g)$.  On $X$, holomorphic sections of
the canonical
bundle are holomorphic  $(2,0)$ forms, and the dimension of the space
$H^{2,0}(X)$ of such forms is the geometric genus $p_g$ of $X$.
We  will always assume that $p_g>0$.
Each   $\alpha\in
H^{2,0}(X)$ can be identified with an element of the $2p_g$-dimensional
real vector space
\begin{equation*}
\H =\ \mbox{Re}\big(\,H^{2,0}\oplus  H^{0,2}\,\big).
\end{equation*}
Using the metric,  each $\a\in \H$ defines an endormorphism
$K_{\alpha}$ of $TX$ by the equation
\begin{equation}
\label{defK}
\langle u,K_{\alpha}v\rangle=\alpha(u,v).
\end{equation}
These endomorphisms $K_\a$ are central to our  discussion,
and we will frequently use the following properties.
Denote by $\nabla$ the Levi-Civita connection of the given  metric.
\begin{lemma}
\label{lemma1}
The $K_\a$ are skew-adjoint and anti-commute with $J$ ($K_\a J=-J K_\a$).
Furthermore,
\begin{equation*}
(a)\ \ \nabla K_\a\ =\ K_{\nabla \a}
\hskip 30pt \mbox{and} \hskip 30pt
(b)\ \ K_\a^2\ =\   - |\a|^2 Id.
\end{equation*}
Consequently, they satisfy the pointwise Clifford relations
$$ K_\a K_\beta+K_\beta K_\a\ =\ -2\langle\a,\beta\rangle\,Id. $$
\end{lemma}
\pf The first two statements and (a) are immediate from
(\ref{defK}). The Clifford relations follow by polarization from
(b), which  is easily  proved (cf.  \cite{l1}). \qed
\bigskip

Now consider holomorphic maps $f:C\to X$ from a connected complex
curve with complex structure $j$ into $X$.  It is standard in geometric
analysis to consider solutions of the perturbed $J$-holomorphic  map equation
$$
\del_Jf=\nu
$$
where $\del_Jf=\frac12(df+Jdfj)$ and where
$\nu$ is an appropriate perturbation term.
In \cite{l1} the first author observed that,
on a K\"{a}hler surface with $p_g>0$,
there is a natural family of such perturbations parameterized by
$\H$. Specifically, we can consider the pairs
$(f,\alpha)$ satisfying
\begin{equation}
\overline{\partial}_Jf= K_\alpha\partial_J f j.
\label{junhoeq}
\end{equation}
This can equally well be viewed as a set of unperturbed holomorphic
map equations for a family of almost  complex structures $\{J_\a\}$
parameterized
by $\H$.  For each $\a \in \H$ the endomorphism $J K_\a$ is skew-adjoint,
so $Id+JK_{\alpha}$ is
injective, and hence invertible. Thus there
is a family of almost complex structures
\begin{equation}
\label{jalpha}
J_{\alpha} = (Id + JK_{\alpha})^{-1}J\,(Id + JK_{\alpha})
\end{equation}
on $X$ parameterized by $\a\in \H$.  A simple computation shows that
(\ref{junhoeq}) is equivalent to the $J_\a$-holomorphic map equation
\begin{equation}
\del_{J_\a}f=0
\label{junhoeq2}
\end{equation}
for maps $f:C\to X$.
Our stucture theorem for GW invariants will emerge from studying the
solutions of this equation for a fixed $\a\in \H$.  Note that while
$\a$ itself is  holomorphic, the corresponding almost complex
structure $J_\a$ need not be integrable.  On the other hand, $J_\a$
is generally not a generic almost complex structure on  $X$, so the
moduli space of $J_\a$ holomorphic maps does not directly define the
GW invariants.

\vskip1cm

\setcounter{equation}{0}
\section{The Localization Lemma and vanishing results}
\label{section3}
\bigskip

The discussion in this section builds on the following simple
principle about Gromov-Witten invariants.
\begin{Vprin}
\label{prin}
If for {\em some} $\w$-tamed almost complex structure $J$,  a class $A\in H_2(X)$ cannot
be represented by a $J$-holomorphic curve of genus $g$, then
$GW_{g,n}(X,A)$ vanishes.
\end{Vprin}
The  proof is straightforward:  if some $GW_{g,n}(X,A)$ were not zero,
we could choose sequences $\{J_n\}$ of generic almost complex
structures converging to $J$
for which there were $J_n$-holomorphic maps representing $A$.  But
then, by the compactness theorem for pseuo-holomorphic maps,   a
subsequence  of those maps would limit to a  $J$-holomorphic
map representing  $A$, contradicting the assumption.
As a simple application, note that for a K\"{a}hler surface $(X, J)$, any
$J$-holomorphic curve represents a (1,1) class, so $GW_{g,n}(X,A)=0$ unless
$A$ is a (1,1) class.  This observation allows us to restrict
attention to (1,1) classes for all our results.
\begin{lemma}(Image Localization Lemma)
\label{ILL}
Fix a K\"{a}hler surface $(X, J)$ with $p_g>0$ and $\a\in \H$.
If $f:C\to X$ is a  $J_\a$-holomorphic map with connected domain that
represents a (1,1) class $A\neq 0$, then $f$    is
in fact a $J$-holomorphic map  whose image $f(C)$   lies in the support
of the zero divisor  $D_\a$ of $\a$.
\end{lemma}

\pf For any $C^1$ map $f:C\to X$ we have the pointwise equality
\begin{equation}
\label{Junho1.3}
\langle \del f,\ K_\a \partial f j\rangle\ d\text{vol} \ =\  f^*\alpha
\end{equation}
(see Proposition 1.3 of \cite{l1}). Integrating over the domain and using
(\ref{junhoeq}) gives
$$
\int_C |\del f|^2\ =\  \int_C \langle \del f,\ K_\a \partial f
j\rangle \ =\  \int_{C}  f^{*}\alpha.
$$
Because $\alpha$ is closed, the last integral  is the homology
pairing $\alpha[A]$.   This  vanishes on the (1,1) class A because
$\a$ is a linear  combination of $(2,0)$ and $(0,2)$ forms.  Thus
$\del f \equiv 0$ on $C$.  Then
using (\ref{junhoeq}), Lemma~\ref{lemma1} and the equality
$|df|^2=|\del f|^2+|\partial f|^2$, we
obtain
\begin{equation}
\label{ILL2}
0\ =\ \int_C |\del f|^2\ =\  \int_C |K_\a \partial f j|^2 \ =\
\int_{C} |\a|^2\,|df|^2.
\end{equation}
Since $A\neq 0$, there is at least one irreducible component of $C$
with $df\not\equiv 0$.  On each such component $C_i$,  $df$ has
finitely many zeros, so (\ref{ILL2}) implies that $f(C_i)$  lies in
the support of $D_\a$.  Each of the remaining components is taken to
a single point by $f$; since $C$ is connected  those  points also lie
in the support of  $D_\a$.
\qed

\bigskip

Lemma \ref{ILL} leads directly to some striking vanishing results for
GW invariants.  For example, K3 and abelian surfaces have trivial
canonical bundle, so admit  (2,0) forms that vanish nowhere.   Lemma
\ref{ILL} and Principle \ref{prin} then give:
\begin{cor}
\label{K3}
For K3  and abelian surfaces, all GW invariants $GW_{g,n}(X,A)$
vanish for $A\neq 0$.
\end{cor}
We also obtain a vanishing result for  the GW numbers
(\ref{GWnumbers}). This follows from the Vanishing Principle and  the
geometric interpretation of the GW numbers.

\begin{cor}
\label{GWga=0}
On a K\"{a}hler surface $X$ with $p_g>0$, any GW invariant
constrained to pass through (generic) points or circles vanishes.
Equivalently, $GW_{g,n}(X,A) (\mu; \ga^1,\dots, \ga^k)=0$ whenever
one of the $\ga^j$ lies in $H^3(X)$ or $H^4(X)$.
\end{cor}
\pf When $PD(\ga^j)$ is a point or 1-dimensional class, we can fix a
geometric representative $\Gamma_j$  disjoint from $D_\a$.   Then,
if the invariant $GW_{g,n}(X,A) (\ga^1,\dots, \ga^k)$ were not zero,
we could find a sequence $\{J_n\}$ of generic almost complex
structures converging to $J_\alpha$  and $J_n$-holomorphic maps
$\{f_n\}$ representing $A$ with $f_n(x_i)\in \Gamma_i$ for all $i$
and $n$. The compactness theorem would then yield a limit
$J_\alpha$-holomorphic map $f$ satisfying $f(x_j)\in \Gamma_j$,
contradicting Lemma \ref{ILL}. \qed

\bigskip

The Image Localization Lemma allows us to localize the GW
invariants for K\"{a}hler surfaces with $p_{g}>0$.  When $X$ is
such a surface and $\a\in \H$, the support of  the zero divisor
$D_\a$ of $\a$ is a union of disjoint topological components
$D_\a^k$. Lemma \ref{ILL}  implies that, for generic $(J,\nu)$
near $(J_\a,0)$, the image of any $(J,\nu)$-holomorphic map with
connected domain  lies in an open neighborhood $U_k$ of one and
only one of the $D_\a^k$.  Thus the compactified moduli space of
$(J,\nu)$-holomorphic maps representing a non-zero class $A$ is a
disjoint union
\be \label{Mdisjointunion}
\bMgn(X,A)\ =\ \coprod
\  \bMgn(U_k,A_k)
\ee
over all $A_k$ with $(\iota_k)_*A_k=A$
under the inclusion $\iota_k:U_k\to X$.  Note that each $U_k$ is
an open symplectic four-manifold with $H_*(U_k)=H_*( D_\a^k)$. As
in Section 1, the image of each $\bMgn(U_k,A_k)$ under the map
(\ref{evmap}) defines a homology class
\be \label{localInvtdef}
GW^{loc}_{g,n}(D_\a^k ,A_k) \in H_*(\ov{\M}_{g,n}\ti D_k^n)
\ee
that we call the {\em local GW invariant of  $D_\a^k$
for the (non-zero) class $A_k$}. These local invariants depend on
the choice of the canonical divisor $D_{\a}$, rather than on the
choice of $\a$ itself. Indeed, if $\b\in\H$ also has zero divisor
$D_{\a}$, then
$\b=c\,\a$ for some constant.  Thus
$J_{\a}$ and $J_{\b}$  are connected by a path
$J_{t}=J_{\a_{t}}$ with $\a_{0}=\a$ and $\a_1=\b$ for which every
$J_{t}$-holomorphic map lies in the  support  of $D_{\a}$.
The standard corbodism argument then shows that the local invariants
$GW^{loc}_{g,n}(D^k_\a,A_k)$ associated with $J_{\a}$ and $J_{\b}$ are
the same.

We remark in passing that the local invariants (\ref{localInvtdef})
can also be regarded as elements of  the homology of the space
$\ov{\M}_{g,n}(D^k_\a,d_k(A))$   of stable maps into  the curve $D^k_\a$
with degree determined by the equation $(\iota_k)_*A_k=d_k(A)[D^k_\a]$.
 From that perspective,  (\ref{localInvtdef}) is the image of the
local invariant under the homology map induced by the evaluation map
(\ref{evmap}) with $X=D^k_\a$.

Pushing (\ref{Mdisjointunion}) forward under the evaluation map
(\ref{evmap}) and passing to homology shows that, for
$A\neq 0$
\be
\label{SumLocalGWI-1}
GW_{g,n}(X,A) \ =\  \sum_{(\iota_k)_*A_k=A} \,GW^{loc}_{g,n}(D^k_\a,A_k).
\ee
for any choice of the canonical divisor $D_{\a}$.  This formula is
the first step toward our structure theorem.  It shows that the GW
invariants can be expressed as a sum of local contributions
associated with the components of a canonical divisor.

\pagebreak

\setcounter{equation}{0}
\section{Local GW invariants}
\label{section4}
\bigskip

The local invariants in the sum (\ref{SumLocalGWI-1})   depend, at
least {\em a priori},  on the local geometry of $J_\a$ around the
components   of the canonical divisor $D_\a$.  In the rest of this paper
we will  write $D_\a=\sum\, m_kD_k$ and assume that
the $D_k$ are smooth and disjoint.   We will show that the local
invariants depend only on discrete data $g,n,d$ and   the multiplicities
$m_k$.  When $D_k$ is
smooth every map with image in $D_k$ represents a multiple $d$ of
$[D_k]$, so we will write the local invariant (\ref{localInvtdef}) as
$$
GW^{loc}_{g,n}(D_k,m_k,d)
$$
or simply $GW^{loc}_{g,n}(D_k,d)$ when  $m_k=1$.
Then, for $A\neq 0$, equation (\ref{SumLocalGWI-1}) reads
\begin{equation}\label{SumLocalGWI}
GW_{g,n}(X,A) \ =\
\sum_{d[D_k]=A} \,GW^{loc}_{g,n}(D_k,m_k,d).
\end{equation}

Using  arguments like those in the previous section, one can also define
local GW invariants
of some open complex surfaces.  Fix a smooth curve $D$ with canonical bundle
$K_{D}$ and    a line bundle $\pi:N\to D$ satisfying $N^{m+1}=K_{D}$.
The total space of $N$ is a complex manifold; from the exact sequence $0\to
\pi^{*}N \to TN \to \pi^{*}TD\to 0$ we see that its canonical bundle
is
\begin{equation}
\label{KNequation}
K_N\, =\,  \wedge^{2}T^{*}N\, =\, \pi^{*}K_{D}\otimes\pi^{*}N^{*}\, =\,
\pi^{*}N^{m+1}\otimes\pi^{*}N^{*}\, =\, \pi^{*}N^m.
\end{equation}
This bundle   has a tautological section $\a$
whose zero divisor is exactly $mD$.  Regarded as
a section of the canonical bundle $K_N$,  $\a$ is
a holomorphic (2,0)-form on $N$.  The argument used to prove
Lemma~\ref{ILL} then shows that
the image of any $J_{\a}$-holomorphic map into $N$  lies in $D$.
On the other hand, an open  neighborhood $U$ of $D\subset N$ is isomorphic
to some open neighborhood $V$ of the zero section $D_{0}$ in the
projectivization
$\P=\P(N\oplus\co_{D})$ by an isomorphism taking $D$ to $D_{0}$.
The pull-back of the K\"{a}hler form on $V$ by that isomorphism gives
an K\"{a}hler form on $U$.
Thus, for any generic $(J,\nu)$ sufficiently close $(J_{\a},0)$
the moduli space $\M_{g,n}(U,d[D])$ can be compactified by standard geometric analysis techniques.
Taking the image as in (\ref{evmap}) yields homology classes
\begin{equation}\label{LI-T}
L_{g,n}(N,m,d)\ \in \ H_*(\ov{\M}_{g,n}\ti D^n)
\end{equation}
that we call the {\em local GW invariants
of $N$   associated with $mD$} for maps representing $d[D]$, $d>0$.
When $m=1$ we will often write (\ref{LI-T}) as simply  $L_{g,n}(N,d)$.
These local invariants depend
on the zero divisor of $\a$  but  not  on $\a$ itself by the
following reasoning.  Let
$\b$ be a section of the canonical bundle $K_N$,
defined on a  neighborhood $U$ of $D\subset N$, such that
the  zero divisor of $\beta$ is $mD$. Then $\b=h\,\a$ for some
holomorphic function $h$ whose  restriction of $h$ to $D$ is non-zero
constant.  Hence, after
shrinking $U$ if necessary,  $J_a$ and $J_\b$ can be connected by a path
$J_{\a_t}$ where the zero divisor of each $\a_t$ on $U$   is $mD$.
As in the previous section, the usual corbordism argument
then shows that
the local invariants associated with $J_\a$ and $J_\b$ are the same.

A similar corbodism argument gives the following fact.

\begin{lemma}\label{def-inv}
If $\{(N_{t},D_{t})\}_{0\leq t\leq 1}$ be a smooth path of line bundles
satisfying $N_{t}^{m+1}=K_{D_{t}}$ then
\begin{equation*}
L_{g,n}(N_{0},m,d)\ =\ L_{g,n}(N_{1},m,d).
\end{equation*}
\end{lemma}
Thus the local invariants (\ref{LI-T}) depend only on the discrete
data $g,n$, and $d$ and the deformation class of the pair $(N,D)$.

\medskip

\begin{ex} \label{Excep-C}
Consider the line bundle $\O(-1)$ on $\P^{1}$.  The
complex structure $J_0$ on the total space of $\O(-1)$ is not of the
form $J_\a$, but nevertheless has the property that any
$J_0$-holomorphic map representing the class $d[\P^{1}]$ has   image
in the zero section in the total space of the bundle
$\O(-1)$.  The argument used above thus applies for $J_0$ as well as
for $J_\a$, showing  that  $J_0$ itself defines the local GW
invariants $L_{g,n}(\O(-1),d)$.
\end{ex}

\bigskip

We can relate the local invariants of $D_k$ defined in
(\ref{localInvtdef}) with the local invariants of its normal bundle
defined in (\ref{LI-T}), as follows.

\begin{lemma} \label{Main-Lemma}
Let $X$ be a K\"{a}hler surface with $p_{g}>0$ and
$D_{\a}=\sum {m_k}D_k$  be the zero divisor of $\a\in \H$.
If $D_k$ is smooth with normal bundle $N_k$  and $D_k\cap
D_\ell=\emptyset$ for all  $\ell \ne k$, then
\begin{equation*}
GW^{loc}_{g,n}(D_k,m_k,d)\ =\ L_{g,n}(N_k,m_k,d).
\end{equation*}
\end{lemma}
\pf  Fix $D=D_k$. By the adjunction formula, the normal bundle $N$ of $D$
satisfies $N^{m+1}=K_{D}$ with $m=m_k$.
Let $Z$ be the  blow-up of $X\times\cx$ along $D\times \{0\}$.  The
projection $X\times\cx\to\cx$ lifts to a map $p: Z\to \cx$ whose
fibers $Z_{\la}=p^{-1}(\lambda)$ are
isomorphic to $X$ for $\la\ne 0$ and whose central fiber $Z_{0}$ is a
singular surface
$X\cup_{D}\P$ where $\P$  is the ruled surface $\P(N\oplus \co_{D})\to D$
defined by fiber projectivization.
The proper transform of $D\times\cx$  is a smooth divisor
$\tilde{D}\subset Z$  disjoint from the proper  transforms
$\tilde{D}_\ell$ of the other $D_\ell\ti \cx$,
and $\a$ gives rise to a section $\tilde{\a}$ of the canonical bundle
$K_Z$ of $Z$ whose zero divisor is
$m\tilde{D}+\sum_{\ell\neq k}m_\ell \tilde{D_\ell} +\P$.
Now fix a tubular neighborhood $U$ of $\tilde{D}$ that is disjoint from
the $\tilde{D}_\ell$.  Let $\kappa$ be the line bundle  of the  divisor
$ m \tilde{D}$, and let  $\beta\in \Gamma(\kappa)$ be a section
with zero divisor $m\tilde{D}$.
For each $\lambda$, the intersection $U_{\la}=U\cap Z_{\la}$
is a tubular neighborhood of  $D_{\la}=\tilde{D}\cap Z_{\la}$.  The restriction
$$
\kappa_{\la}\ =\ \kappa |_{U_\la}.
$$
is  the line bundle on $U_\la$ with divisor $m D_\la$.  Observe that:
\begin{itemize}
\item For $\la\neq 0$, the normal bundle $N_\la$ to $Z_\la$ in $Z$ is
trivial.  Restricting the exact sequence
$0\to TZ_\la \to TZ\to N_\la \to 0$ to $U_\la$ then shows
that the canonical bundle of $U_\la$ is the restriction of the canonical bundle
of $Z$, which is the bundle of the divisor
$m\tilde{D}\cap U_\la= m D_\la$.
\item  For $\la=0$ we use a different argument.  By the definition of blow-up,
$U_0$ is biholomorphic to a neighborhood of the zero section in the total space
of the bundle $N\to D$; in fact, this identifies the zero section with $D_0$.
But by (\ref{KNequation}) the canonical bundle of  $N$ has a tautological
section whose divisor is $m$ times that zero section.
\end{itemize}
Thus  $\kappa_\la$ is the canonical bundle of $U_\la$ for each $\la$.

\smallskip
Restricting $\beta$ to $U_\la$ gives a section $\beta_\la$ of $\kappa_\la$
whose zero divisor is $m D_{\la}$ and a  corresponding almost complex
structure $J_{\la}=J_{\b_{\la}}$  on $U_{\la}$.  Then
the image of any $J_{\la}$-holomorphic map lies in $D_{\la}$,
so  $J_{\la}$  determines  local invariants $GW^{loc}_{g,n}(D_\la,m,d)$
of $U_{\la}$ for the class $d[D_{\la}]$ (with $d>0$).   Because
$\beta_\la$ and $J_{\la}$ vary  smoothly in $\la$, we then have
\begin{equation*}
GW^{loc}_{g,n}(D_\la,m,d)\ =\ GW^{loc}_{g,n}(D_0,m,d)
\end{equation*}
for each $\la$. The righthand side of the above equals
$L_{g,n}(N,m,d)$ by definition, while for $\la\neq 0$ the lefthand
side is $GW^{loc}_{g,n}(D,m,d)$ because $Z_\la$ is biholomorphic to $X$
by a map that takes $D_\la$ to $D$.
This completes the proof of the lemma.
\qed

\medskip

\begin{ex} Let $\pi:E(m+2)\to \P^1$ be an elliptic surface with
$12(m+2)$ singular fibers which are all nodal. This surface is K3 if
$m=0$ and properly elliptic if $m>0$.
By the canonical divisor formula\,(see (\ref{eq1}) below) the canonical
bundle of $E(m+2)$ is $\pi^*\O(m)$.  Thus
the generic  canonical divisor is the sum of $m$ disjoint regular fibers $F_i$,
and  for any regular fiber $F$ the divisor $mF$ is also a canonical
divisor. Using  Lemma~\ref{Main-Lemma}  and equation (\ref{SumLocalGWI})
we then have
$$ L_{g,n}(\O,m,d[F])  \ =\ GW_{g,n}(E(m+2),d[F])\ =\  m\,L_{g,n}(\O,d[F]). $$
\end{ex}

\vskip 1cm

\setcounter{equation}{0}
\section{Exceptional curves and blowups}
\label{section5}
\bigskip

This section establishes a  ``blowup formula''  that reduces
the problem of computing GW invariants
to the case of  minimal surfaces.
This extends some previous partial blowup formulas, cited at the end
of this section.
In our approach the blowup formula is a consequence of the
localization Lemma~\ref{ILL}.

First consider a closed  symplectic 4-manifold   $X$ with an  almost
complex structure $J$  and an exceptional $J$-holomorphic curve $E$.
We can then consider the (global) invariants
\bear
\label{localE}
GW_{g,n}(X,d[E])
\eear
which gives the contributions to $GW_X$ of all maps whose image
represents a multiple  of $[E]$.
Fix a diffeomorphism $\iota:\P^1\to E$ and let $\iota_*$  denote the
map $H_*(\bMgn\ti (\P^1)^n)\to H_*(\bMgn\ti X^n)$ induced by $\iota$.
\begin{lemma}
\label{exceptionalcurvelemma}
For $d>0$, (\ref{localE}) is given by the local invariant of
Example~\ref{Excep-C}:
\begin{equation*}
GW_{g,n}(X,d[E])\ =\ \iota_*L_{g,n}(\co(-1),d).
\end{equation*}
\end{lemma}
\pf  Since $E^2=-1$, any $J$-holomorphic curve representing a class
$d[E]$ has image in $E$.  Thus
$$
GW_{g,n}( X,d[E])\ =\  GW_{g,n}^{loc}( E,d).
$$
After rescaling the symplectic form on $\P^1$ we may assume that
$\iota:\P^1\to E$ a symplectomorphism.  By the Symplectic
Neighborhood Theorem this extends to a symplectomorphism
$\varphi:U\to V$ from a neighborhood $U$ of  the zero section in $\O(-1)\to
\P^1$ to a neighborhood $V$ of $E$ in $X$.
Pushing the standard complex structure $J_0$ on $\O(-1)$ forward by
$\varphi$  gives an almost complex structure $J'_0$  on $V$ that
makes $\varphi$ an isomorphism of almost complex neighborhoods.
Furthermore,  $E$ is a $J_0'$ holomorphic curve,  so the local
invariant above can be calculated using $J_0'$.  Thus when $d>0$
\begin{equation*}
GW_{g,n}^{loc}(E,d)\ =\ \iota_*L_{g,n}(\co(-1),d).
\mbox{\qed}
\end{equation*}

\bigskip

Let $X$ be a compact K\"{a}hler surface with $p_g>0$ and let
$\pi:\tilde{X}\to X$ be the blowup of $X$ at a point $p$.  Different
choices of the point $p$ yield surfaces $\tilde{X}$ that are
symplectic deformation equivalent, so the
GW invariants of $\tilde{X}$ are independent of the choice of $p$.
Note that
every $A\in H_2(\tilde{X})$ can be uniquely written as $A= B+dE$
where $E$ is the class of the exceptional curve and $B\cdot E=0$
and the invariant $GW_{g,n}(X,\pi_{*}B)$ can be regarded as
a homology class in $H_{*}(\CM_{g,n}\times (X\setminus \{p\})^{n})$.

\bigskip

\begin{prop}
\label{blowupProp}
Let $X$ be a compact K\"{a}hler surface with $p_g>0$ and let
$\pi:\tilde{X}\to X$ be its blowup at a point $p$.
Then  the  GW invariant of each class $A= B+dE$  as above
is given by
\be
\label{BlowupFormula}
GW_{g,n}(\tilde{X}, A)=
\left\{
\begin{array}{ll}
L_{g,n}(\O(-1), d)  & \mbox{if $A=dE$ with $d>0$}\\[.1cm]
\pi_*^{\prime}GW_{g,n}(X, \pi_*A) & \mbox{if $A\cdot E=0$}\\[.1cm]
0 & \mbox{otherwise}
\end{array}
\right.
\ee
where $\pi_*^{\prime}$ is the induced homology map by the composition
of the isomorphism
$X\setminus\{p\}\to \tilde{X}\setminus E$ and the inclusion
$\tilde{X}\setminus E\to \tilde{X}$.
\end{prop}
\pf  Fix a holomorphic $(2,0)$ form $\alpha$ on $\tilde{X}$ with zero
divisor  $D\in |K_X|$ and a blowup point $p\notin D$.  Then
$\tilde{\alpha}=\pi^*\alpha$ is a holomorphic $(2,0)$ form on
$\tilde{X}$ whose zero divisor  $\tilde{D}\in |K_{\tilde{X}}|$ is the
disjoint union of $\pi^{*}(D)$ and the exceptional curve $E$.
Each class $A\in H_2(\tilde{X})$ with non-zero GW invariant can be
represented by a $J_{\tilde{\alpha}}$-holomorphic map $f:C\to
\tilde{X}$ from a connected curve $C$.  By Lemma~\ref{ILL},
the image of $f$ lies in $\tilde{D}$. Hence either $A\cdot E=0$ or $A=dE$ with
$d>0$.  The case  $A=dE$ was done in Lemma~\ref{exceptionalcurvelemma}.

If $A\cdot E=0$,  choose a sequence of almost complex
structures $J_\ell$  converging to $J_{\tilde{\alpha}}$.  As $\ell\to
\infty$, the $J_\ell$-holomorphic maps converge pointwise to
$J_{\tilde{\alpha}}$-holomorphic maps.  These limit maps lie in
$\tilde{D}$ but not in $E$ because of the condition $A\cdot E=0$.
Thus for large $\ell$ the images are bounded away from $E$; in fact, they are
uniformly bounded away from $E$ for $f$ in  the compact  space
$\bM^{J_\ell}_{g,n}(\tilde{X},A)$ of stable maps.  Consequently, the
condition that $J_\ell$ is generic for this space of stable curves is
the same as the condition that an almost complex structure that
agrees with $\pi_*J_\ell$ outside a sufficiently small neighborhood
of the blowup point is generic for the corresponding space of stable
maps into $X$.  When both are generic,  composition with $\pi$ gives
a diffeomorphism
$$
\M^{J_\ell}_{g,n}(\tilde{X},A) \overset{\approx}{\to
}\M^{J_\ell}_{g,n}(X,\pi_*A)
$$
that respects orientations and  the stablization and evaluation maps.
Hence the
corresponding GW invariants are equal.
\qed

\bigskip

\begin{rem}
The  hypothesis $p_g>0$ is needed in Proposition \ref{blowupProp} .
For example, when $X$ is $\P^2$ and $L$ is  the class of the line,
the invariants  $GW_{g,n}(\tilde{X}, aL+bE)$ with $b>1$ are non-zero:  they are
enumerative counts of the curves in $\P^2$ satisfying certain contact
and tangency
conditions at the blowup point (see Gathmann \cite{Gath}).
Jianxun Hu showed  that the part of  Proposition~\ref{blowupProp} pertaining to classes $A$ with $A \cdot E= 0$ and $A
\cdot E= 1$
hold on any symplectic manifold (\cite{Hu}).  For other classes, however, the
contrast between Proposition \ref{blowupProp}  and Gathmann's results for
$\P^2$ shows that any universal blowup formula for GW invariants must
distinguish  rational surfaces from those with $p_g>0$.
\end{rem}
The first part of the Structure Theorem~\ref{introthm}  is a version
of the blowup formula (\ref{BlowupFormula}).  Given a compact \Ka
surface $X$ with $p_g>0$,  let $\pi:X\to X'$ be the projection to the
minimal model.   By perturbing the  blowup points, we can insure that
there is a canonical divisor on $X$  whose support is a disjoint
union of exceptional curves $\{E_k\}$ and other curves $D_\ell$.
Define a formal power series with coefficients in $H_*(\bMgn\ti
(\P^1)^n)$ by setting
\bear
\label{defLE}
L^0(t)\ =\
\sum_{d>0}\sum_{g,n}\,\frac{1}{n!}\,L_{g,n}(\O(-1),d)\,t^{d}\,
\la^{2g-2}
\eear
and another with coefficients in $H_*(\bMgn\ti X^n)$ by
\bear
\label{defGW'}
GW^{\prime}_{X}\ =\ \sum_{A\ne 0}\sum_{g,n}\, \frac{1}{n!}\,
\pi^{\prime}_{*}GW_{g,n}(X,\pi_*A) \,t_A\, \la^{2g-2}.
\eear
The blowup formula then gives the following succinct equation (cf.
Theorem~\ref{introthm}).

\begin{prop}
\label{exceptionalcurveProp}
The  GW invariant of $X$ is a sum
$$
GW_X\ =\ GW^0_X\ +\ \sum_{E_i} L^0(t_{E_{i}})\ +\ GW'_X.
$$
\end{prop}
\vskip1cm

\setcounter{equation}{0}
\section{The Structure Theorem for properly elliptic surfaces}
\label{section6}
\bigskip

In light of the blowup formula of  the previous section, we can
henceforth assume that all surfaces $X$ are minimal.   Furthermore,
the GW invariants of  a K3 or abelian surface are trivial  by
Corollary~\ref{K3}.
The  Enriques-Kodaira
classification then shows that, among minimal  surfaces with $p_g>0$,
there are two cases left to consider:   minimal properly elliptic
surfaces and minimal surfaces of general type.  We will consider
these separately.

\medskip

Let $\pi:X\to C$ be a minimal properly elliptic surface. Then the
sheaf $L=(R^1\pi_*\O_X)^{-1}$ is a line bundle on $C$ with $\deg
L=\chi(\O_X) \geq 0$, and the canonical bundle is
$$
K_X\ =\ \pi^*(L\otimes K_C)\ \otimes\ \O(\sum_k\,(m_k-1)F_k^\prime)
$$
where  $F_k^\prime$ are multiple fibers of multiplicity $m_k$
(\cite{FM} pages 47-49).  Correspondingly, each canonical divisor of
$X$ has the form
\begin{equation}\label{eq1}
  \ \sum_j\,n_jF_j \ +\ \sum_k\,(m_k-1)F_k^\prime
\end{equation}
where $\sum n_jF_j$ is the pullback of a divisor in $|L+ K_C|$
of degree $k_\pi=\chi(\O_X)-2\chi(\O_C)$.  In general, the
  fibers $F_j$ need not to be  smooth or disjoint from the
$F_k^\prime$.

\begin{prop}
\label{FriedmanProp}
Every minimal properly elliptic surface $\pi:X\to C$ can be deformed
  to a minimal properly  elliptic surface whose generic canonical
divisor has the form (\ref{eq1}) where  the $F_j$ and $F_k^\prime$
are disjoint  smooth fibers.
\end{prop}

\pf  By a theorem of Moishezon  $X$ can be deformed to a minimal
properly elliptic surface whose only singular fibers are reduced
nodal curves and multiple fibers with smooth reduction (see \cite{FM}
page 113 and \cite{BHPV} page 266).
This deformed $X$ is obtained by  log transforms on an elliptic surface
$\pi:S\to C$ without multiple fibers  whose canonical bundle is  $K_S\ =\ \pi^*(L+K_C)$ for
the same line bundle $L$ (\cite{FM} pages 102-103).  By deforming the fibers
on which the logarithmic
transformations are done we can assume that none of the fibers $F_k^\prime$
lie over the base points of
the linear system $|L+K_C|$, and hence   the generic canonical divisor of
$X$ has the form (\ref{eq1})
with $F_j\cap F_k^\prime=\emptyset$ for all $j$ and $k$.
 It therefore suffices
to prove the  Proposition~\ref{FriedmanProp} for the surface $S$.

Next note that  $|L+K_C|$ is empty when $\deg
(L+K_C)=\deg L+2g(C)-2 \leq 0$ and has  a base point at $p\in C$  if
and only if $h^0(L+K_C-p)=h^0(L+K_C)$ (see \cite{H} page 308). By
Riemann-Roch and Serre duality, this last condition is equivalent to
$ h^0(p-L)=h^0(-L)+1$.   Hence
$|L+K_C|$ has no base points   when $\deg L\geq 2$, and also when
$\deg L=1$ and $L \neq\O(p)$ for any $p\in C$.  In these cases
Bertini's Theorem implies the generic canonical divisor is the
disjoint union of  distinct  smooth fibers.
This leaves only  two specific cases:
\begin{enumerate}
\item[a)]  $\deg L=0$ and $g=g(C)\geq 2$, and
\item[b)]  $L=\O(p)$ for some $p\in C$ and $g\geq 1$.
\end{enumerate}
In fact,
case a) occurs only when  $S$ has no singular fibers (\cite{FM}, page
48).  Thus the proposition is true in case a).

In case b),  choose points $p, q\in C$ that
are not linearly equivalent,
and let $L$ be any
one of the $2^{2g}$  line bundles on $C$ with $L^2=\O(p+q)$.  Following
\cite{FM} page 60, one can construct an elliptic surface
$\pi_L:S_L\to C$ with section with $(R^1\pi_{L*}\O_{S_L})^{-1}=L$
whose only singular fibers are the fibers over $p$ and $q$.  It follows from
Seiler's Theorem (Corollary I.5.14 of \cite{FM}) that each $S_L$ is
deformation
equivalent to $S$.   Since $L$ is not isomorphic to $\O(p)$ or $\O(q)$,
  the generic element of  $|L+K_C|$ has support disjoint from $p$ and
$q$.  The corresponding canonical divisor of $S_L$ is then
 a union of smooth fibers.
\qed

\begin{remark}{\rm
R. Friedman (private communication) has proved a stronger version
of Proposition~\ref{FriedmanProp}: one can assume,  after further
deformations, that  $n_j=1$ for all $j$.   This is a
more natural statement, but is not needed for our purposes in light
of the calculation of Example~4.4.}
\end{remark}

 Proposition~\ref{FriedmanProp} is useful because
 {\em  K\"{a}hler surfaces that are deformation equivalent as complex  surfaces have the same GW
invariants}.
 This is true because deformation equivalent surfaces are smoothly isotopic
 (\cite{FM} page 18) and, because the space of K\"{a}hler  forms
 with a fixed orientation is convex,
that isotopy lifts to give a symplectic deformation equivalence.  Consequently, the
 GW invariants are the same.

\bigskip

Thus we may assume that the
generic canonical divisor $D$ has the form (\ref{eq1}) where
\begin{itemize}
\item each $F_j$ is a regular  fiber with
holomorphically trivial normal bundle, and
\item smooth  multiple fiber $F_{m_k}$ whose normal bundle $N_{k}$ is
torsion of order
$m_{k}$ in the group
$\mbox{Pic}^{0}(F_{m_{k}})$ of line bundles of degree zero
\end{itemize}
(cf. Section III.8 of \cite{BHPV}).
Then for a regular fiber  $F_j$ with $n_j=1$,
we have the local GW invariants
(\ref{LI-T}) with $m=1$ and $N=\O$.  These define a function $L^1$ as
follows.

\begin{defn}
\label{DefLF}
Let $\O$ is the trivial line bundle over $T^{2}$ and set
\begin{equation*}
L^{1}(t)\ =\
\sum_{d>0}\sum_{g,n}\,\frac{1}{n!}\,L_{g,n}(\O,d)\,t^{d}\,\la^{2g-2}.
\end{equation*}
\end{defn}
For a regular fiber  $F_j$ with $n_j>1$
one can form the corresponding power series  with  $L_{g,n}(\O,d)$
replaced by $L_{g,n}(\O,n_j,d)$.
The result is simply $n_jL^1(t)$ by the calculation of Example~4.4.

\vskip.6cm

For multiple fibers, we will define similar functions $L^2_m(t)$
in terms of the GW invariants of a   ``model space'' constructed by a
logarithmic
transformation. To that end, fix
an elliptic $K3$ surface $X\to\P^1$, a regular $F$ of $X$ and a
torsion line bundle $\xi\in Pic^0(F)$ of order $m>1$.  Applying the
logarithmic transformation defined by this data
yields an elliptic surface $X(F,\xi)$.  This surface
\begin{itemize}
\item   is simply connected  and therefore K\"{a}hler\,(see \cite{go}
and Theorem 3.1 of \cite{BHPV}), and
\item  has $\chi(\O_X)=2$, so  by (\ref{eq1}) its
canonical divisor  $D= (m-1)F^\prime_m$
is supported on a single multiple fiber $F_m^\prime$ of multiplicity $m$.
\end{itemize}
Changing the choices of $X$, $F$ and $\xi$ yields a surface that is
deformation equivalent to $X(F,\xi)$  (Theorem I.7.6 of \cite{FM}) and hence has the same GW
invariants.  We will write $K3(m)$ for the generic surface in this
deformation class.

\begin{defn} \label{6.defL3}
With $K3(m)$ and $F_m^\prime$ as above,  set
\begin{equation*}
L^{2}_m(t) \  = \
\sum_{d>0}\sum_{g,n}\,\frac{1}{n!}\,GW_{g,n}(K3(m),d[F^\prime_m])\,t^d\,\la^{2g-2}.
\end{equation*}
\end{defn}

The following proposition shows that the local invariants
at any any smooth multiple
fiber  $F_m$   of multiplicity $m$
can be expressed in terms of GW invariants of
$K3(m)$ that are encoded in  the function $L^2_m(t)$.

\begin{prop}
\label{FmTheorem}
Let $X$ be a properly elliptic surface with a smooth multiple
fiber $F_m$ of multiplicity $m\geq 2$.
Then $$ GW_{g,n}^{loc}(X,m-1,d[F_m])\ =\
   GW_{g,n}(K3(m),d[F_m^\prime]) $$
\end{prop}

\pf Recall that  there is a local model
for a neighborhood $U$ of $F_m$  (cf. Prop. 6.2 of [FM]).
Specifically, there is a (smooth) elliptic
fibration $\pi_0:U_0\to \Delta$ over a unit disk $\Delta\subset \cx$
and a torsion line bundle
$\xi$ of order $m$ on $\pi_0^{-1}(0)$
such that $U$ is isomorphic, as an elliptic fibration, to
the elliptic fibration obtained by performing  the $m$-logarithmic
transformation defined by $\xi$ on the central fiber $\pi_0^{-1}(0)$.
In particular, $\pi_0$ and  $\xi$ completely
determine the $m$-spin
curve $(F_m, N_m)$, that is, determine  the curve $F_m$  and a normal bundle $N_m$ satisfying $N_m^{m+1}=K_{F_m}$.

Furthermore, there is a holomorphic function $h_0$ on $\Delta$
satisfying $\mbox{Im}\, h_0(s) >0$ such that
$\pi_0:U_0\to \Delta$ is the quotient
$(\cx\times \Delta)/(\Z\times \Z)\to \Delta$
with the action of $\Z\times \Z$  given by
$$ (m,n)(z,s)\ =\ (z+m+nh_0(s),s) $$
(pg. 202 of \cite{BHPV}).
Now fix a normal neighborhood of a smooth fiber of $K3\to \P^1$.
One can then  choose an isomorphic (smooth) elliptic fibration
$\pi_1:U_1\to \Delta$ over the unit disk $\Delta$
under which the fixed smooth fiber of $K3$ corresponds to
the central fiber $\pi_1^{-1}(0)$. As above,
this fibration is determined by a holomorphic function $h_1$ on $\Delta$
with $\mbox{Im}\, h_1(s) >0$.

Since for each $t\in [0,1]$ the function $h_t=(1-t)h_0+th_1$ is holomorphic on
$\Delta$ and satisfies $\mbox{Im}\, h_t(s) >0$, using $h_t$
one can obtain a family of elliptic fibrations $\pi_t:U_t\to \Delta$.
Then, performing $m$-logarithmic transformation on each fiber
$\pi_t^{-1}(0)$ using a family of line bundles $\xi_t$ of order $m$ on
$\pi_t^{-1}(0)$ with $\xi_0=\xi$ shows that the $m$-spin curves defined by
$F_m \subset X$ and a  multiple fiber $F_m^{\prime}\subset K3(m)$ are deformation equivalent. Therefore, we have
\begin{equation*}
GW^{loc}_{g,n}(X,m-1,d[F_m])\ =\ GW^{loc}_{g,n}(K3(m),m-1,d[F^{\prime}_m]) \ =\ 
   GW_{g,n}(K3(m),d[F^{\prime}_m])
\end{equation*}
where the first equality follows from Lemmas 4.1 and 4.3, and
the second follows from (3.5) and the fact that the canonical divisor
of $K3(m)$ is $(m-1)F^\prime_m$.
\qed

\bigskip

The structure  theorem for minimal properly elliptic surfaces follows
immediately from  (\ref{SumLocalGWI}),  Lemmas~\ref{def-inv} and
\ref{Main-Lemma},  and Proposition~\ref{FmTheorem}.  The
result is the following case of Theorem~\ref{introthm}.

\begin{theorem}\label{E-str}
If  $X$ is a minimal properly elliptic surface whose canonical
divisor  $D$ is given as in (\ref{eq1}),
then
$$
GW_X\ =\ GW^0_X\ +\ k_\pi L^{1}(t_F)\ +\   \sum_k
L^{2}_{m_k}(t_{F_k})
$$
where $F$ is a regular fiber and $t^{m_k}_{F_k}=t_F$.
\end{theorem}

\vskip1cm

\setcounter{equation}{0}
\section{The Structure Theorem for surfaces of general type}
\label{section7}
\bigskip

When $X$ is a minimal surface of general  type, every canonical
divisor is connected and has arithmetic genus $h\geq 2$
(\cite{BHPV}).   Unlike the case of  elliptic surfaces,  it is not
always possible to deform a surface of general  type to insure the
existence of a smooth  canonical divisor. For example, Bauer and
Catanese have shown that there is a surface $S$ with
$p_g = 4$, $K^2 = 45$ that has  no complex deformations and such that
each canonical divisor is singular and reducible (\cite{BC}).  It is not presently understood how common such examples are. To avoid this complication we make the  following assumption.

\bigskip

\noindent{\bf Assumption.}\ \ {\em For some \Ka  structure in the
deformation class of  $X$, there is a smooth  canonical  divisor $D$
with multiplicity 1.}

\medskip
\noindent  (Of course, if this is true for some \Ka structure then it
is true for the generic one.)  When $D$ is smooth with  multiplicity 1,  the adjunction formula
shows that the normal bundle $N$ of $D$ is a holomorphic square root
of  $K_D$:
\bear
\label{adjunction}
N^2= K_{D}.
\eear

\bigskip

Recall that a  {\em theta characteristic} on a smooth curve $D$ is a line
bundle $N$ with $N^{2}=K_{D}$.  In the special case when $K_D=\O$ is trivial,
the set $S(D)$ of all theta characteristics is the same as the group
$J_{2}(D)$ of points of order 2 in the Jacobian.  In general, $S(D)$
is a principal homogeneous space for $J_{2}(D)$ with the obvious
action:  if $N$ is a theta characteristic and $L^2=\O$ then $N\otimes
L$ is another theta characteristic.  Since $J_{2}(D)$ is naturally
isomorphic to $ H^{1}(D;\Z_{2})$, there are  $2^{2h}$ theta characteristics on a curve of genus $h$.
A theta characteristic $N$ is {\em even} or {\em odd} according to
the parity of
$h^{0}(D,N)$.

A {\em spin curve} is a pair $(D,N)$ consisting of a curve with a theta characteristic.
 The spaces ${\mathcal S}_{h,+}$ (resp. ${\mathcal S}_{h,-}$) of  all genus $h$  even (resp. odd) spin curves have compactifications
$\overline{\mathcal S}_{h,\pm}$ The following  three   facts are
classical.

\begin{prop}\label{T:parity}
(see \cite{ACGH}, \cite{A} and \cite{C})\  Let $D$ be a smooth curve of genus $h$.
\begin{enumerate}
\item[(a)]  There are $2^{h-1}(2^h+1)$ even and $2^{h-1}(2^h-1)$ odd
             theta characteristics.
\item[(b)]  $h^0(D_{t},N_{t})\ \mbox{mod 2}$ is constant along any smooth
             family $(D_t,N_t)$ of spin curves.
\item[(c)] $\overline{\mathcal S}_{h,\pm}$ is an irreducible projective variety and $\partial {\mathcal S}_{h,\pm} = \overline{\mathcal S}_{h,\pm}\setminus {\mathcal S}_{h,\pm}$ is a proper analytic subvariety.
\end{enumerate}
\end{prop}

\medskip

\begin{cor}
\label{NewSpinCor}
The invariants $L_{h,n}(N,d)$, defined by (\ref{LI-T}) when $m=1$, depend only on the genus $h$ and the parity of $h^{0}(D,N)$.
\end{cor}
\pf Since $\overline{\mathcal S}_{h,\pm}$ is irreducible, the smooth
part $\overline{\mathcal S}_{h,\pm}^*$ is connected (\cite{GH} page
21), and hence $\overline{\mathcal S}_{h,\pm}^*\setminus
\partial{\mathcal S}_{h,\pm}$ is connected.  Thus any two smooth spin
curves of the same parity can be joined by a path of spin curves.
The Corollary then follows from Lemma~\ref{def-inv}. \qed

\bigskip

For our case --- a surface of general type with a smooth  canonical  divisor $D$
with multiplicity 1 --- the parity of $h^{0}(D,N)$ is actually a {\em global} invariant, as the following lemma shows.

\medskip

\begin{lemma}\label{chi}
If $X$ is a  minimal surface of general type and
$D\subset X$ is  a smooth canonical  divisor  with normal bundle $N$, then
$$
h^{0}(D,N)\equiv\chi(\co_{X}) \ {\rm (mod \  2)}.
$$
\end{lemma}
\pf
Since $N$ is the restriction of $K$ to $D$, there is an exact sequence
$0\to {\mathcal O}_{X} \overset{m}{\to} {\mathcal O}_{X}(K)
\overset{r}{\to} {\mathcal O}_{D}(N) \to 0$,
where $m(f)=f \a$ and $r(\b)=\b_{|_{C}}$. This  induces a long exact
sequence of cohomology which, using the isomorphisms
$H^{0,1}(X)\cong H^{1}({\mathcal O}_{X}) $ and $ H^{1}(K)\cong
H^{2,1}(X)$, begins
\begin{equation*}
0 \to\ H^{0}({\mathcal O}_{X})\ \to\ H^{0}(K)\ \to\ H^{0}(N)\ \to\
H^{0,1}(X)\ \overset{m_*}{\to} \ H^{2,1}(X)\ \to\ \cdots
\end{equation*}
where  $m_*$ is given by $m_*(\la)=\la\wedge \a$. The hermitian inner
product on $H^{0,1}(X)$ gives an orthogonal splitting
$H^{0,1}(X)=\ker m_* \oplus V$ and, by the above sequence, $h^{0}(N)
=  p_{g}+q-1-\dim V$.  Since $\chi(\co_{X})=1-q+p_g$,  it suffices to
show that $V$ is
even dimensional.
After composing with the star operator, $L=*m_*:
H^{0,1}(X)\to H^{0,1}(X)$ satisfies
\begin{equation}
\label{Lequation}
\langle\, \la\,,L(\delta)\,\rangle\ =\
-\langle\, \delta\,,L(\la)\,\rangle.
\end{equation}
Thus $L$ induces a nondegenerate sympletic pairing  on $H^{0,1}(X)/\ker m_* \cong V$, so $\dim V$ is even.
\qed

\bigskip

We  can proceed as we did for elliptic
surfaces. Again, we first define invariants associated with a spin curve.
\begin{defn}
\label{defL}
For each smooth genus $h\geq 2$ spin curve $(D,N)$ and let $L_{g,n}(N,d)$ be the local GW invariant
(\ref{LI-T}) and set
\bear
\label{6.gtype}
L^3_{h,\pm}(t)\ =\ \sum_{g,n}\sum_{d\geq 1}
\frac{1}{n!}\,L_{g,n}(N,d)\ t^d\lambda^{2g-2}
\eear
This notation incorporates the fact that, by Corollary~\ref{NewSpinCor},  this series depends only on $h$ and the parity of   $(D,N)$.
\end{defn}
For minimal surfaces of general type, the statement of the structure
theorem is especially simple because the canonical divisor of $X$ has
a single component.  The GW series is obtained from one of the series
(\ref{6.gtype}).

\begin{theorem}
Suppose that $X$ is a minimal surface of general type with a smooth,
multiplicity 1 canonical  divisor $D$.  Let $h=K_X^2+1$ be the genus
of $D$. Then (again suppressing inclusion maps)
\begin{equation*}
GW_{X}\ =\  GW_{X}^{0}\ +\
\begin{cases}
 L^3_{h,+}(t_{D'}) & \mbox{if $\chi(\O_X)$ is even}\\
  L^3_{h,-}(t_{D'}) & \text{if $\chi(\O_X)$ is odd.}
  \end{cases}
\end{equation*}
Consequently, the GW series of $X$ depends only on $h$ and $\chi(\O_X)$.
\end{theorem}

\pf This follows directly from (\ref{SumLocalGWI}), Lemma~\ref{Main-Lemma}, Corollary~\ref{NewSpinCor}
and Lemma~\ref{chi}.
\qed

\vskip 1cm

\setcounter{equation}{0}
\section{Moduli spaces and  linearizations}
\label{section8}
\bigskip

For each fixed $\a\in\H$, we can consider the linearization $D_f$ of
the $J_\a$-holomorphic map equation at each
$J_\a$-holomorphic map $f:C\to X$.  This operator is important for
local descriptions of the moduli space.  After a brief discussion of
moduli spaces, we will write down the formula for $D_f$ and show that it has
some remarkable
analytic properties.

Consider a smooth component $D$ of a canonical divisor of $X$.  When $D$  has  multiplicity 1,  we have $N^2= K_{D}$ as in (\ref{adjunction}).
When $D=F_m$ is a multiple elliptic fiber with multiplicity $m$, the
normal bundle satisfies $N^{m}=\O_D$.  Taking chern classes,  both cases give
the formula
\bear
\label{adjunction2}
c_1(N)[D]=h-1.
\eear
\begin{lemma}
\label{dimformula}
Fix a smooth  genus $h$ component $D\subset X$ of a canonical
divisor.  Then the (formal real) dimensions of the moduli spaces
$\M_g(D,d)$ (of degree $d$ genus $g$ covers of the curve $D$)  and
$\M_g(X,d[D])$ (of maps from a genus $g$ curve representing $d[D]\in
H_2(X)$) are
\be
\label{dimM24}
\dim\M_g(D,d)\ =\ 4\beta
\qquad\mbox{and}\qquad
\dim\M_g(X,d[D])\ =\ 2\beta
\ee
where $\beta=d(1-h)+g-1$.
\end{lemma}
\pf  The restriction of $TX$ to $D$  decomposes as $TD \oplus N$.
Using (\ref{adjunction2}) we then have $K_X\cdot D = K_D\cdot D
-c_1(N)[D] = h-1$. Both parts of (\ref{dimM24}) then follow from the
dimension formula (\ref{dimM}).
\qed

\bigskip

To interpret the number $\beta$ geometrically, consider a $J_\a$-holomorphic map $f:C\to D$ from a smooth genus $g$ curve onto  $D$.  The canonical classes of $C$ and $D$ are then related by
the Riemann-Hurwitz formula $K_C =f^*K_{D} +B$ where $B$ is the
ramification divisor. Consequently, the number of branch points,
counted with multiplicity, is
\bear
\label{RHeq}
|B| =2\beta\qquad\mbox{where}\ \beta=d(1-h)+g-1.
\eear

\bigskip
To proceed, we need explicit formulas.
By a standard calculation (cf. \cite{IS}, \cite{rt1}),
the   linearization   of the $J_\a$-holomorphic map equation, evaluated at a map $f$ and
applied to a variation $\xi$ of the map and a variation $k$ of
the complex structure on the domain, is
\bear
\label{fulllinearization}
D_f(\xi, k)\ =\ L_f(\xi) + J_\a df k
\eear
where  the
operator $L_{f}: \Omega^0(f^*TX) \to \Omega^{0,1}(f^*TX)$  given by
$$
L_{f}(\xi)(w)\, = \, \del_f \xi (w)+
\left(\frac12J\nabla_\xi J+\nabla_\xi K_\alpha\right)( df jw)
+K_\alpha(\nabla\xi)jw
$$
for each $w\in\Omega^0(TC)$ (here  $\del_f\xi(w)$ is $
\frac12\left(\nabla_w \xi+J\nabla_{jw}\xi\right)$).
In our case $\nabla J=0$ and  $\alpha$ vanishes along the image of
$f$, so that
\bear
L_f= \del_f+R_\alpha
\label{linearization}
\eear
with
\bear
R_\alpha(\xi)\ =\  (\nabla_\xi K_\alpha)\circ df\circ j.
\label{defRa}
\eear

\begin{lemma}\label{lemmaX}
Let $D$ be a smooth component of a canonical divisor $D_\a$ and
$N$ be the normal bundle of $D$. Then,
for each $p\in D$, $u\in T_p\,D$ and $\xi\in N_p$ we have
\begin{center}
\begin{tabular}{ll}
(a)\  $\nabla_u K_\a =0$ \hspace*{5.4cm}& (c)\ $\nabla_{J\xi} K_\a (u)
= -J\nabla_\xi K_\a (u)$\\ \\
(b)\  $\nabla_\xi K_\a(u)$ is orthogonal to $T_p\,D$ &
(d)\  $|\nabla_\xi K_\a (u)|^2 = |\nabla \a|^{2}|\xi|^2|u|^2$\\ \\
\end{tabular}
\end{center}
\end{lemma}
\pf Since $\a\equiv 0$ along $D$, (a) follows from Lemma~\ref{lemma1}a.
Next,  the fact $\a$ is a closed 2-form and $\nabla_u\a=0$ gives the formula
\begin{equation*}
0\ =\ d\alpha(u,\xi,\eta)\  =\  (\nabla_\xi\alpha)(\eta,u) \ -\
(\nabla_\eta\alpha)(\xi,u)
\end{equation*}
for any $\eta\in T_p\,X$.  Applying the definition of $K_\a$, this becomes
\begin{equation}\label{E:Ka}
\langle\eta,\nabla_\xi K_\alpha (u)\rangle\ =\
\langle\xi,\nabla_\eta K_\alpha (u)\rangle .
\end{equation}
When $\eta\in T_p\,D$, we have
$\nabla_\eta K_\a=0$ and thus (\ref{E:Ka})
shows (b).
Then, because  $J$ is skew-adjoint with $\nabla J=0$,
and $K_\alpha J=-J K_\alpha$, (\ref{E:Ka}) implies that
\begin{equation*}
\langle\eta,\nabla_{J\xi} K_\alpha (u)\rangle
= \langle J\xi,\nabla_\eta K_\alpha u\rangle
= \langle \xi,\nabla_\eta K_\alpha (Ju)\rangle 
= \langle\eta,\nabla_\xi K_\alpha (Ju)\rangle =
-\langle\eta,J\nabla_\xi K_\alpha (u)\rangle.
\end{equation*}
This gives (c). Finally, noting $\nabla K_\a(Ju)=-J\nabla K_a(u)$
and using (c), the fact $TD$ and $N$ are $J$-invariant and
Lemma~\ref{lemma1}, we have
\begin{equation*}
|\nabla_\xi K_\a(u)|^2\ =\ |\nabla K_\a|^2|\xi|^{2}|u|^2\ =\
|\nabla \a|^2|\xi|^{2}|u|^2.
\mbox{\qed}
\end{equation*}

\pagebreak

As an immediate corollary, we have\,:
\begin{cor}
\label{Rlemma1}
If  $f:C\to D$ is a $J_\a$-holomorphic map
onto a smooth component $D$ of a canonical divisor $D_\a$,
then $R_{\a}$ vanishes on
$f^*TD$ and defines a (real) bundle map
$$
R_\a:f^*N\to T^{0,1}C\otimes f^*N
$$
where $N$ is the normal bundle to
$D$.  This $R_\alpha$  satisfies $R_\a J= - J R_\a$ and
\begin{equation}
\label{rasquared}
|R_\a(\xi)|^2\ =\  |\nabla\a|^2\, |\xi|^2\, |df|^2.
\end{equation}
\end{cor}

\bigskip

On a \Ka surface $X$, each $\alpha\in\H$ has an associated
almost complex structure $J_\a$ and canonical divisor $D_\a$.
Let  $V$ is a smooth component of the support of $D_\a$.
Following \cite{ip1}, one can use the space of $(J_\a,\nu)$-holomorphic maps
to define the relative GW invariant for the pair $(X,V)$ provided $(J_\a, \nu)$
is a generic ``$V$-compatible'' pair as defined in Section 3 of \cite{ip1}.

\begin{cor}
\label{VcompatCor}
Let  $V$ is a smooth component of the support of $D_\a$.
If $\nabla\a \equiv 0$ on $V$ then
$(J_a,0)$ is a $V$-compatible pair, while if $\nabla\a \not\equiv 0$ then
$(J_\a, \nu)$ is  $V$-compatible for no choice of $\nu$.
\end{cor}
\pf
Let $\pi_N$ denote the orthogonal projection onto the normal bundle $N$ of $V$.
A pair  $(J_\a, \nu)$ is  $V$-compatible if it satisfies three conditions\,:
$J_\a$ preserves $TV$, $\nabla J_\a$  satisfies
\begin{equation}\label{con-b}
\pi_N\left[(\nabla_\xi J_\a + J_\a\nabla_{J_\a\xi} J_\a)(u)\right]\ =\
\pi_N\left[(\nabla_u J_\a + J_\a\nabla_{J_\a u} J_\a)(\xi)\right]
\end{equation}
for all $u\in TD$ and $\xi\in N$,
and $\nu$ and $\nabla \nu$ satisfy conditions that are automatically true
when $\nu=0$. Since $\a= 0$ along $V$, the definition (\ref{jalpha})
of $J_\a$ shows that
$J_\a=J$ and $\nabla J_\a=-2\nabla K_\a$ at each point in $V$.  Thus
$V$ is $J_{\a}$-holomorphic.
One can then use Lemma~\ref{lemmaX} to see that Condition
(\ref{con-b}) is equivalent to
\begin{equation*}
\nabla_{\xi}K_{\a}(u)\ =\ 0\ \ \forall u\in TD,\  \forall\xi\in N.
\end{equation*}
Lemma~\ref{lemmaX}\,d then implies that $V$-compatibility conditions hold
only if $\nabla\a=0$ along  $V$, and that if   $\nabla\a=0$ along  $V$ then
$(J_a,0)$ satisfies the $V$-compatibility conditions.
\qed

\vspace{.5cm}

The two terms of the operator (\ref{linearization}) satisfy a
remarkable property under the $L^2$ pairing:

\begin{lemma}
\label{Radjjoint}
Let $D$ be a smooth component of a canonical divisor $D_\a$ with
normal bundle $N$.  Then for each
$J_\a$-holomorphic map $f:C\to D$ we have
$$
\int_C \langle \del\xi,\ R_\a\eta\rangle\ +\ \int_C \langle \del
\eta,\ R_\a\xi\rangle\ =\ 0
\quad \ \mbox{$\forall\, \xi, \eta\in \Omega^0(f^*N)$}.
$$
\end{lemma}
\pf Let $f_{s,t}$ be a 2-parameter family of deformations of the map
$f=f_{0,0}$ with  $\frac{d\ }{ds}f|_{s=t=0}=\xi$ and $\frac{d\
}{dt}f|_{s=t=0}=\eta$.  Then  $\del f_{s,t} =\del(s\xi+t\eta)
+Q(s,t)$ where $Q$ is at least quadratic in $(s,t)$.   Since the
image of $f$ represents a multiple of the (1,1) class $[D]$, equation
(\ref{Junho1.3}) gives
$$
0\ =\ \int f^*\a\ =\ \int \langle \del f,\ K_\a(df j)\rangle
$$
for each $f=f_{s,t}$.  Now differentiate this equation  with respect to
both $s$ and $t$ and evaluate at $s=t=0$, noting that $\del f$ and
$K_\a(df j)$ both vanish at $s=t=0$.  The result is
$$
0\ =\ \int \langle \del\xi,\ \nabla_\eta K_\a(df j)\rangle\ +\ \int
\langle \del \eta,\ \nabla_\xi K_\a(df j)\rangle.
$$
The lemma follows by the definition of $R_\a$.
\qed

\bigskip

We finish this section by discussing the operator given by  the normal  component of the linearization (\ref{linearization}).  For each map $f:C\to D$  as in Corollary~\ref{Rlemma1}  the pullback  $f^*TX$ of the tangent bundle decomposes orthogonally as $ f^*TX=f^*TD\oplus f^*N$.   Let $\pi^N$ be the projection onto $f^*N$.   The normal component $\pi^N\circ\nabla$ of the connection on $f^*TX$ is a hermitian connection on $f^*N$; its (0,1) part defines an operator $\del_f^N$ and hence a holomorphic structure on $f^*N$.  The  restriction of $\del_f$ to $f^*N$ then has  the form
$$
\left.{\del_f}\right|_{f^*N} =  \del^N_f +A
$$
where $A$ is a bundle map $f^*N\to T^{0,1}C\otimes f^*TD$ (which vanishes if $f^*N$ is a holomorphic subbundle; see \cite{GH} pg. 78).  On the other hand, since
$f^*TD$ is a holomorphic subbundle,  the restriction of $\del_f$
to $f^*TD$ is an operator $\del_f^T$ on $f^*TD$ which is the usual $\del$-operator.
Corollary~\ref{Rlemma1} then implies that the linearization
(\ref{fulllinearization}), as an operator
$$
D_f : \Omega^0(f^*TD\oplus f^*N)\oplus H^{0,1}(TC) \to
\Omega^{0,1}(f^*TD\oplus f^*N),
$$
is given by
\begin{equation}\label{fulllinearization2}
D_f =\begin{pmatrix}
\del_f^T & A\\0&L^N_f
\end{pmatrix} \oplus Jdf
\end{equation}
where $L^N_f=\del_f^N+R_\a$.
The next result shows that $L^N_f$ is injective.

\medskip
\begin{prop}
\label{DelRinvertible}
Suppose that $f:C\to D$ is a $J_\a$-holomorphic map from a smooth
curve onto a smooth component $D$ of   a canonical divisor $D_\a$ and
that either (i) $\nabla \a \neq 0$ somewhere on $D$, or (ii) $\ker\,
\del^N_f = 0$.     Then
\begin{equation}
\label{kerLf}
{\rm ker}\,  L^N_f= 0.
\end{equation}
\end{prop}
\pf Suppose there is a non-zero  $\xi\in{\rm ker}\,  L^N_f$.  Then the integral
$$
\| L^N_f\xi\|^2\ =\
\|(\del^N_f+R_\alpha)\xi\|^2\ =\
\int_C  |\del^N_f\xi|^2\ +\ |\nabla\a|^2\, |\xi|^2\, |df|^2
$$
vanishes (here we have used  (\ref{rasquared}) and noted that, because $R_\alpha\xi$ is normal and $A\xi$ is tangent,  Lemma
\ref{Radjjoint} holds with $\del_f\xi=\del_f^N\xi+A\xi$ replaced by $\del_f^N\xi$).  But both $\xi$ and $f$ satisfy elliptic equations, so by the Unique
Continuation Theorem for elliptic equations $|\xi|^2\,
|df|^2$ is not zero on any open set.   We conclude that  $\del^N_f\xi=0$ and $\nabla\a\equiv 0$ along $D$.
\qed

\vskip 1cm
\setcounter{equation}{0}
\section{Zero-dimensional spaces of stable maps}
\label{section9}
\bigskip

The simplest GW invariants are those associated with a space of
stable maps whose formal dimension is zero.   Such stable maps  are
especially simple: Lemma~\ref{lemma8.1} below shows that they are
unramified maps from smooth domains, and that  the linearization
$D_f$  is invertible.   Thus all zero-dimensional
GW invariants are signed  counts of the number of connected etale covers.
This section establishes some basic facts needed to make these counts. Specific computations are done in
Section~\ref{section10}.

The  formal dimension of a space $\Mgn(X,A)$ of stable maps is the
index of linearization $D_f$ at each $f\in\Mgn(X,A)$. Calculating as
in the proof of Lemma~\ref{dimformula}
one finds that $\ind D_f=2\beta+2n$ and similarly $\ind L_f^N= -2\beta+2n$.
Consequently, when the space of stable maps is formally 0-dimensional and the domain curve is smooth, we have
\bear
\label{9.1}
\ind D_f=\ind L^N_f=0.
\eear

 Now fix  a smooth canonical divisor satisfying the conditions of Prop\-osition~\ref{DelRinvertible}.  The Image Localization Lemma~\ref{ILL} implies that all invariants $GW_{g,n}(X,A)$ vanish unless $A$ is a multiple $d[D]$ of the class of a component $D$ of that canonical divisor.
 These invariants also vanish  whenever the formal dimension of $\bM_{g,0}(X,A)$ is negative because the space  $\bM_{g,0}(X,A)$, and therefore $\bM_{g,n}(X,A)$,  is then empty for generic $(J,\nu)$.
Thus, using the dimension formula of Lemma~\ref{dimformula}, we may assume that $A=d[D]$ and
\bear
\label{9.2}
\beta=n=0
\qquad\mbox{with}\quad
\beta=d(1-h)+g-1.
\eear
where $h$ is the genus of $D$.

\begin{lemma}
\label{lemma8.1}
Suppose that $D\subset X$ is a smooth component of a canonical
divisor $D_\a$.    Then any non-constant stable map $f:C\to D$ satisfying (\ref{9.2})  is an etale cover  from a smooth curve $C$ and the linearization $D_f$ is invertible.
\end{lemma}
\pf   By (\ref{9.2}) we have $g=dh-d+1$.
Suppose that $C$ has $\ell$ irreducible components $\{C_i\}$.
Restricting $f$ to each component and lifting to the normalization
gives maps $\tilde{f}_i:\tilde{C}_i\to D$.  Suppose that exactly $k$
of these have degree $[\tilde{f}_i]=d_{i}>0$. Then $\sum\,d_{i}=d$
and  (\ref{RHeq}), applied to each $C_i$, gives
\begin{equation*}
g\, =\, dh-d+1\, \leq\, \sum\,(d_{i}h - d_{i} + 1 + \b_{i})\,
=\, \sum\,g_{i}^{\prime} \, \leq\ \sum\,g_{i}\, \leq\, g
\end{equation*}
where $\b_{i}$ is the ramification index of  $\tilde{f_{i}}$,
$g^{\prime}_{i}$ is the geometric genus of $C_{i}$,
and $g_{i}$ is the arithmetic genus of $C_{i}$.
This shows that $k=1$ and $C_{1}$ has the same
geometric and arithmetic genus. Consequently, $C_1$ is smooth
of genus $g$ and the remaining $\ell-k$ components
have genus 0. Stability then implies that $\ell-k=0$.
Thus $C$ is smooth  and $f:C\to D$  has no critical points.

Recall that the linearization $D_f$ is given by (\ref{fulllinearization2}).
The normal operator $L^N_f$ is injective by
(\ref{kerLf}) and hence is surjective by (\ref{9.1}).  Furthermore, $Jdf$ induces an isomorphism from
$ H^{0,1}(TC)$ to  $H^{0,1}(f^*TD)=\cok \del_f^T$, and therefore
$$
\del_f^T\oplus Jdf: \Omega^0(f^*TD)\oplus H^{0,1}(TC) \to \Omega^{0,1}(f^*TD)
$$
is also onto.  Thus $D_f$ is surjective with index zero, so is an isomorphism between the appropriate Sobolev spaces.
\qed

\bigskip

When  $D_f$ is invertible, there is an associated invariant:  its mod
2 spectral flow.
That spectral flow is computed in the next proposition.
This calculation is crucial to the discussion in the next section.

The mod 2 spectral flow of $D_f$ is determined by choosing a path
$D_t$  of first  order elliptic operators from an invertible complex
linear operator $D_0$ to $D_1=D_f$ so that $D_t$ is invertible except
at finitely many $t_i$ along the way, and taking
\
\be
\label{SFdef}
SF(D_f)=\sum_i \dim \ker D_{t_i}\qquad \mbox{(mod 2)}.
\ee
\
This is a homotopy invariant of the path, and is independent of $D_0$
because any two choices of $D_0$ can be connected by a path $D_t$ of
complex linear first  order elliptic operators, and at each point
along such a path $\ker D_t$ is even-dimensional.

\medskip

\begin{prop}
\label{relativesign3}
Under the conditions of Lemma~\ref{lemma8.1},
\be
\label{signLf}
SF(D_f)\equiv  h^0(f^*N)\quad \mbox{(mod 2)}.
\ee
\end{prop}

\pf  First deform $D_f$ to a diagonal operator along the path
$$
D_t =\begin{pmatrix}
\del_f^T & tA\\0&L^N_f
\end{pmatrix} \oplus Jdf.
$$
Because both  $\del_f^T\oplus Jdf$ and $L^N_f$ are surjective, each $D_t$ is surjective with index zero, so $\ker D_t=0$ for all $t$.  Noting that $\del_f^T\oplus Jdf$ is  complex-linear,  we then have
$$
SF(D_f)=SF(D_0)=SF(L_f^N).
$$
Next, since $L_f^N$ is invertible by Lemma \ref{DelRinvertible},
$SF(L_f^N)=SF(L_f^N+B)$ for any sufficiently small compact perturbation $B$.
Now write $L^N_f$ as $\del +R_{\a}$ with  $\del= \del^N_f$.
Because $\ind \del =0$ for etale covers, we can choose a
complex-linear isomorphism $\bar{B}:\ker \del \to \cok\del$ and set
$B=\bar{B}P$ where $P$ is the $L^2$ orthogonal projection onto $\ker
\del$.  Then
$$
D_t= \del +\delta B +tR_\a
$$
is a path from $D_0=\del +\delta B$ to $D_1=L^N_f+\delta B$.  Using
Lemma \ref{Radjjoint}, we have
$$
\int |D_t\xi|^2\ =\ \int |\del \xi|^2 +|(\delta B+tR_\a)\xi|^2.
$$
This shows that $D_0$ is invertible and that  $\ker D_t$  lies in
$\ker \del$ and in  $\ker (\delta \bar{B}+tR_\a$) for each $t$.
Taking  $\delta$ sufficiently small, we then have
$$
SF(L_f^N)\ =\ SF(D_1)\ =\ SF(\delta \bar{B}+\bar{R_\a})
$$
where $\bar{R_\a}$ is the restriction of $R_\a$ to $\ker\del$.   But
$\bar{R_\a}$ is injective  and anti-commutes with $J$ by  Lemma
\ref{Rlemma1}. Furthermore,  its image is $L^2$ perpendicular to the
image of $\del$ by Lemma \ref{Radjjoint} and  $\ind \del =0$, so
$\bar{R_\a}:\ker \del \to \cok\del$ is an isomorphism.  This means
that $SF(\delta \bar{B}+\bar{R_\a})$ is the same as $SF(\bar{R_\a})$
and, from the definition (\ref{SFdef}), the same as
$SF(\bar{B}^{-1}\bar{R_\a})$.  Here $\bar{B}^{-1}\bar{R_\a}$ is an
isomorphism of $H^0(C,f^*N)$ that anti-commutes with $J$.  The lemma
is completed using   two simple  facts about the spectral flow of
finite-dimensional matrices:
\begin{enumerate}
\item[(a)] $(-1)^{SF(A)} = \sign \det A$ for all $A\in GL(n,\R)$.
\item[(b)]  If $A\in GL(2n,\R)$ satisfies $JA=-AJ$ then $SF(A)=n$ mod 2.
\end{enumerate}
To see (a), choose a path $A_t$ in the space of $n\times n$ matrices
from $A$ to $Id$; for a generic such path each kernel in
(\ref{SFdef}) is 1-dimensional, so the spectral flow is the number of
sign changes in $\det A_t$.  For (b), choose a basis $\{v_1, Jv_1,
\dots, v_n, Jv_n\}$ and set $w_i=Av_i$.  Then $v_1\wedge
Jv_1\wedge\dots \wedge J v_n$ and
$w_1\wedge Jw_1\wedge\dots \wedge J w_n$  both represent the complex
orientation, so the calculation
\begin{align*}
\det A\cdot v_1\wedge JV_1\wedge\dots \wedge J v_n \ &=\ Av_1\wedge
AJv_1\wedge\dots \wedge AJ v_n\\
&=\ (-1)^n w_1\wedge Jw_1\wedge\dots
\wedge J w_n
\end{align*}
shows that $\sign \det A=(-1)^n$.
\qed

\bigskip

In Gromov-Witten theory,   the GW invariant  associated  with a
zero-dimensional space of stable maps is the signed count of the maps
in that space with the sign of each map $f$  specified by the mod 2
spectral flow of the linearization $D_f$ (provided each $D_f$ is an
isomorphism). By Proposition \ref{relativesign3} this sign is
\bear
\label{etalesignformula}
(-1)^{SF(D_f)}\ =\  (-1)^{h^0(f^*N)}.
\eear

This sign is well-defined even though $h^0(f^*N)$ may change under deformations of  the holomorphic structure on
$f^*N$.   This is because, for etale covers $f:C\to D$,  we have $f^*K_D=K_C$ and hence  the equation
$N^2=K_{D_0}$ pulls back to $(f^*N)^2=K_C$.   Thus $(C, f^*N)$ is a spin curve, so by Lemma~\ref{T:parity} the parity of $h^0(f^*N)$ does not change as $(D,N)$ is deformed.

\medskip

Formula (\ref{etalesignformula}) is a key difference between GW invariants in two and four dimensions.   The finite set of  etale covers of $D$ contribute to both  the Gromov-Witten invariants of the curve $D$, and to the GW invariants of $X$ through the inclusion $D\subset X$.  But in the first case each etale cover contributes $+1/|\mbox{Aut($f$)}|$ to the invariant, while in the second case the signs vary according to (\ref{etalesignformula}).

\vskip 1cm
\setcounter{equation}{0}
\section{Zero-dimensional GW invariants:  computations}
\label{section10}
\bigskip

The facts  established in the previous
section are enough to compute the contributions of
etale covers to the GW series in some cases.  We do this  for the
canonical class itself, for double covers, and for general etale
covers for elliptic fibers.

\pagebreak
\noindent{\bf{The canonical class}}

\medskip
When $X$ and $D$ are as in Lemma~\ref{chi}, $D$ is an
embedded  genus $g=K^{2}+1$  curve  representing the canonical class
$K$. For that genus, the GW invariant has dimension 0 by (\ref{dimM})
and  is immediately computable using (\ref{SumLocalGWI}),
Proposition~\ref{relativesign3}, and
Lemma~\ref{chi}\,:
\begin{equation*}
GW_{g}(X,K)\ =\ GW^{loc}_{g}(D,1)\ =\ (-1)^{h^{0}(N)}\ =\
(-1)^{\chi(\co_{X})}.
\end{equation*}
This fact is well-known  from other perspectives.  In the context of
Taubes' $Gr$ invariant (see \cite{T}), $g=K^{2}+1$ is the ``embedded
genus'' case.  In that case the $Gr$ invariant is the same as the
Seiberg-Witten invariant and is  given by
$Gr(K)=SW(K)=(-1)^{\chi(\co_{X})}$. On the other hand,  because $D$
is embedded and connected,
we also have $Gr(K)=GW_{g}(X,K)$.

\vskip.5cm
\noindent{\bf{Double covers}}

\medskip

The etale double covers of a curve $D$ are classified by either
$H^{1}(D;\Z_{2})$ or equivalently by $J_{2}(D)$.  In fact,  if the
square of a line bundle $L$ is trivial, then $L$ has a bisection $s$
satisfying $s^2=1$ and the image of $s$ is a smooth unramified
double covering $f:C_{L}\to D$ that is connected whenever
$L\neq\O_D$. Such double coverings satisfy $f_{*}{\mathcal
O}_{C_{L}}={\mathcal O}_{D}\oplus L^{-1}$ and thus for any line
bundle $N$ on $D$
\begin{align}\label{push}
h^{0}(C_{L},f^{*}N)\ &=\ h^{0}(D,f_{*}f^{*}N)\ =\ h^{0}(D,N\otimes
f_{*}{\mathcal O}_{\tilde{C}_{L}}) \\
&=\ h^{0}(D,N)\ +\
h^{0}(D,NL^{-1}). \notag
\end{align}
Now suppose that $D$ is a smooth component of a canonical divisor of
genus $h$ with normal bundle satisfying $N^{2}=K_{D}$. Since each map
$f$ in the moduli space $\CM(D,2)$ of etale double covers with
connected domains has automorphism group $\Z_2$,
each contributes $\pm \frac12$  to the GW invariant, with the sign
given by Proposition~\ref{relativesign3}.     Thus
Proposition~\ref{T:parity}a, Lemma~\ref{chi} and equation
(\ref{push}) yield
\begin{equation}
\label{deg=2}
GW_{g}^{loc}(D,2)\  =\  
\sum_{f\in \CM(D,2)}\frac{1}{2}\,(-1)^{h^{0}(f^{*}N)} 
\  = \   \frac{1}{2}\left[(-1)^{h^{0}(N)}\,2^{h}\ -\ 1\right]
\end{equation}
where $g=2h-1$.   For surfaces of general
type the sign $(-1)^{h^0(N)}$ can be calculated from the global invariant
$\chi(\co_{X})=1-q+p_g$  by Lemma~\ref{chi}.

\begin{example}
\label{DoubleCoverExs}{\rm
Exceptional curves have no etale double covers, while elliptic fibers
have three connected  double covers, all  etale with genus 1.  Thus
(\ref{deg=2}) gives
\begin{enumerate}
\item  A regular fiber $F$ has trivial normal bundle, so
$GW^{loc}_{1}(F,2) = -\frac32.$
\item A multiple fiber $F_2$ of order 2 has $h^0(N)=0$, so  $
GW^{loc}_{1}(F_2, 2) =\frac12$.
\item Formula (\ref{deg=2}) does not apply to multiple fiber $F_m$
with multiplicity $m>2$ because the normal bundle to $F_m$ is not a
theta characteristic,
but instead satisfies $N^{m}=\O$.  Nevertheless,   we have $h^{0}(f^{*}N)=0$
for each of the three nontrivial double covers $f$ of $ F_{m}$, so
$$ GW^{loc}_{1}(F_m,m-1,  2)= \frac32. $$
\item When $D$ is  a smooth  multiplicity 1 canonical divisor in  a
surface of general type, $D$ has genus $h=K^{2}+1$ and a
connected double cover $C\to D$ is etale if and only if  $C$ has
genus $g=2K^{2}+1$.
By (\ref{deg=2}) the genus $g=2h-1$ invariant invariant is
$$
GW_{g}(X,2K)\ =\    GW_{g}^{loc}(D,2)\ =\ \frac{1}{2}\left[
(-1)^{\chi(\co_{X})}\,2^{h}\ -\ 1\right].
$$
\end{enumerate}
}
\end{example}

\vskip.4cm
\noindent{\bf{Etale Covers of  Elliptic Fibers}}

\medskip

When $(X,J)$ is a generic complex structure on a minimal properly
elliptic surface, the generic  canonical divisor has components  of two
types:  smooth elliptic fibers and multiple fibers with smooth reduction.    The simplest cases are regular fibers and multiple fibers of multiplicity two.  For those, we can give  explicit formulas for the contributions to the GW invariants of smooth etale covers.

\bigskip

\noindent{\em \bf Regular Fibers.}\ \    Every holomorphic map $f:C\to
F$ onto a  regular elliptic fiber has $f^*N=\O$, so $h^0(f^*N)=1$.
Such a map $f$  is an etale cover  if and only if $C$ has genus $g=1$.
The stable moduli space $\bM_{1,0}(F,d)$  consists of
$\sigma(d)$ points, where  $\sigma(d)=\sum_{k|d} k$ is the sum of
the  divisors of $d$.  Each of these is a generic as  $J_\a$-holomorphic map
(Lemma~\ref{lemma8.1} implies that $\ker D_f=0$)
with  automorphism group of order $d$,  and each  is counted with a minus sign
by Lemma \ref{relativesign3} because $f^*N=O_C$.
Thus the contribution of the  etale covers to the local GW invariant
of  $F$ is
$$
\sum_{d>0} GW^{loc}_1(F,d)\ t_F^d\ =\
-\sum_{d>0}\frac{\sigma(d)}{d}\ t_F^d\ =\ -\int
\frac{G(t_F)}{t_F}\,dt_F
$$
where
$$
G(t)\ =\   \sum\,\sigma(d)\,t^d \ =\ \prod_{k>0}
\frac{k t^k}{1-t^k} .
$$

\bigskip

\noindent{\em \bf $F_2$ Fibers.}\ \  As in (\ref{eq1}), every elliptic fiber $F_2$ of multiplicity 2 is a
component  of each canonical divisor $D_\a$ with multiplicity 1.  In particular,  $\nabla\a$ does
not vanish  identically along $F_2$.
Thus  by Lemma~\ref{lemma8.1}  $\ker D_f=0$  and $\cok D_f=0$  for
every $J_{\a}$-holomorphic
etale cover $f:C\to F_2$.
Consequently, the (local) GW invariants of  etale covers are determined by
their Taubes' type. Since the degree 1 map has positive sign,  and
two of (nontrivial)
double covers have positive sign and one has negative sign, we have
\begin{align*}
\sum_{d>0} GW^{loc}_1(F_2,d)\ t_{F_2}^d\ &=\ \sum_{d>0}\,
\frac{1}{d}\,\left[ \sigma(d)\ -
2\,\sigma\Big(\frac{d}{2}\Big)\right]\,t^{d}_{F_{2}} \\ &=\
\int \frac{G(t_{F_{2}})-2G(t^{2}_{F_{2}})}{t_{F_2}}\,dt_{F_2}
\end{align*}
(see Proposition 4.4 of \cite{ip0}).

\vskip 1cm
{\small

\vskip1cm

\noindent {\em Department of Mathematics, Michigan State University,
East Lansing, MI 48824}

\smallskip

\noindent {\em E-mail addresses:}\ \ {\ttfamily leejunho\@@msu.edu}\\
\phantom{E-mail addresses:\ } {\ttfamily parker\@@math.msu.edu}

}
\end{document}